\documentclass[a4]{article}
\usepackage{amsmath}

\font\es=eufm10

\def\gC{\mbox{\es {C}}}

\def\ov{\overline}

\def\dfrac#1#2{\displaystyle \frac{#1}{#2}}

\def\H{\mbox{\boldmath $H$}}

\def\N{\mbox{\boldmath $N$}}

\def\R{\mbox{\boldmath $R$}}

\def\0{\mbox{\boldmath {0}}}    
\def\1{\mbox{\boldmath {1}}}      
\def\2{\mbox{\boldmath {2}}}      
\def\3{\mbox{\boldmath {3}}}      
\def\4{\mbox{\boldmath {4}}}      
\def\5{\mbox{\boldmath {5}}}      
\def\6{\mbox{\boldmath {6}}}      
\def\7{\mbox{\boldmath {7}}}      
\def\8{\mbox{\boldmath {8}}}      
\def\9{\mbox{\boldmath {9}}}      
\def\a{\mbox{\boldmath $a$}}

\title{On the conjugacy of pure imaginary elements of quaternion algebras and Cayley algebras }
\author{Takashi \textsc{Miyasaka}}
\date{}
\begin{document}
\maketitle
\begin{abstract}
In the algebras  $\H$, $\H'$, $\H^{C}$, $\gC$, $\gC^\prime$ and $\gC^C$, we show some results on the conjugacy of two pure imaginary non-zero elements with same norm.
\end{abstract}
\begin{flushleft}
\large{\bf Introduction }
\end{flushleft}

It is known that any automorphism $\alpha$ of the field $\H$ of quaternions is an inner automorphism and any two pure imaginary quaternions $a$, $b$ with the same norm are conjugate, that is, there exists a  pure imaginary non-zero element $p  \in \H$ such that $b=pap^{-1}$ (Theorem 6.(1)). In the case of the split quaternion algebra $\H'$ (resp. the complex quaternion algebra $\H^{C}$), any two pure imaginary non-zero elements $a,b$ with the same norm are conjugate, that is, there exists an invertible element $p \in \H'$ (resp. $\H^{C}$) such that $b=pap^{-1}$ (Theorem 6.(2)). 
In the case of the division Cayley algebra $\gC$, any two pure imaginary elements $a,b$ with the same norm are conjugate, that is, there exists a pure imaginary non-zero element $p \in \gC$ such that $b=pap^{-1}$ (Theorem 3.(1)).
In the case of the split Cayley algebra $\gC^\prime$ (resp. the complex Cayley algebra $\gC^C$), we need twice conjugate operations, that is, for any two pure imaginary non-zero elements $a,b$ with the same norm, there exist invertible pure  imaginary elements $p,q \in \gC^\prime$ (resp. $\gC^C$) such that $b=q(pap^{-1})q^{-1}$ (Theorem 3.(2)).
We obtain the above results  in the constructive manner, by using concrete elements. 

Finally, the author would like to thank Takae Sato, Osamu Shukuzawa and Ichiro Yokota for their earnest guidance, useful advice and  constant encouragement.
\vspace{4mm}

\begin{flushleft}
\large{\bf 0. Preliminaries }
\end{flushleft}

Let $\H$ be the field of quaternions with the canonical $\R$-basis $\{1,\,e_1,\,e_2 ,\,e_3 \}$ with the usual multiplication defined by 
$$
\begin{array}{c}
e_1{}^2 =e_2{}^2 =e_3{}^2 = -1,\\
 e_1e_2 = e_3=-e_2e_1, \quad e_2e_3 = e_1=-e_3e_2, \quad e_3e_1 = e_2=-e_1e_3, 
 \end{array}
 $$
and let $\H^C$ be the complexfication algebra of $\H$: $\H^C=\H \oplus i\H$. The algebra $\H'$ of split quaternions is defined as follows. $\H'$ is the algebra with $\R$-basis $\{1,\,e_1{}',\,e_2 ,\,e_3{}' \}$  with the multiplication defined by 
 $$
 \begin{array}{c}
 e_1{}'^{\,2} =e_3{}'^{\,2} =1,e_2{}^2 = -1,\\
 e_1{}'e_2 = e_3{}'=-e_2e_1{}', \quad e_2e_3 {}'= e_1{}'=-e_3{}'e_2, \quad e_3{}'e_1{}' = -e_2=-e_1{}'e_3{}'.
\end{array} 
 $$
Next, let $\gC = \H \oplus \H e_4$ (resp. $\gC' = \H' \oplus \H' e_4 $) be the division Cayley algebra (resp. the split Cayley algebra) over $\R$ with the multiplication
$$
\begin{array}{c}
   (m_1 + n_1e_4)(m_2 + n_2e_4) = (m_1m_2 - \ov{n_2}n_1) + (n_1\ov{m_2} + n_2m_1)e_4,
 \end{array}
 $$ 
 for $m_1 + n_1e_4, m_2 + n_2e_4 \in \gC$ (resp. $\gC'$), where  $\overline{m}$ is the conjugate element of $m \in \H$ (resp. $\H'$), 
and let $e_5=e_1e_4,e_6=-e_2e_4,e_7=e_3e_4$ (resp. $e_5{}'=e_1{}'e_4,e_6=-e_2e_4,e_7{}'=e_3{}'e_4$).
The complex Cayley algebra $\gC^C$ is defined as the complexification of $\gC$: $\gC^C=\gC  \oplus i \gC$. In $\gC$, $\gC'$ and $\gC^C$, the conjugation is defined by $\overline{m+ne_4}=\overline{m}-ne_4$. 
In the algebras $K=\H,\, \H',\, \H^C,\, \gC \, \,, \,\gC'$ and $\gC^C$ above,  the inner product $(a,b)$ and the norm $\N(a)$ by 
$$
\begin{array}{l}
(a, b) = \dfrac{1}{2}(a\overline{b} + b\overline{a}), \quad \N(a)=(a,a)=a \overline{a}.
\end{array}
 $$ 
Note that if $a \in K$ satisfies $\N(a) \neq 0 $, then $a$ is invertible and the inverse element $a^{-1}$ of $a$ is given by $a^{-1}=\overline{a}/{\N(a)}$. Finally, we use the following notation
$$
\begin{array}{l}
K_0=\{a \in K \mid \overline{a}=-a \}, \quad \quad \,\, K_0{}^*=\{a \in K_0 \mid a \neq 0 \}, 
\vspace{1mm}\\
K{}^{\times}=\{a \in K  \mid \N(a) \neq 0 \},  \quad K_0{}^{\times}=\{a \in  K_0 \mid \N(a) \neq 0 \}.
\end{array}
$$
For $a \in K_0$, $\N(a)$ is nothing but $-a^2$.
 
 \vspace{3mm}
\begin{flushleft}
{\large{\bf 1. Cases of Cayley algebras}}
\end{flushleft}
\vspace{3mm}


{\bf Lemma 1.} (1) {\it For any $a \in (\gC)_0$, there exists $p \in {(\gC)_0}^*$ such that $pap^{-1}=-a $.}\\
\quad (2) {\it Let $K =\gC', \gC^C$. For any $a \in K_0$, there exists $p \in {K_0}^\times$ such that $pap^{-1}=-a $.}
\vspace{2mm}

{\bf Proof.} We may assume that $a \neq 0$.\\
\quad (1) Case $a \in  (\gC)_0$. Express $a = a_1e_1 + \cdots + a_7e_7, a_k \in \R$. Then, at least one element $p$ of the following elements 
$$
\begin{array}{l}
   a_2e_1 - a_1e_2, \quad a_3e_2 - a_2e_3, \quad a_5e_4 - a_4e_5,
	\quad a_7e_6 - a_6e_7
\end{array}
$$
satisfies $pap^{-1}=-a$. Indeed, each element $p$ above satisfies $pa=-ap$, and the norm $\N(p)$ are
$$
\begin{array}{l}
  {a_2}^2+{a_1}^2, \quad {a_3}^2+{a_2}^2,\quad {a_5}^2+{a_4}^2,\quad {a_7}^2+{a_6}^2.
\end{array}
$$
If $\N(p)=0$ for all $p$, then we have $a_1=\cdots=a_7=0$, that is, $a=0$, which contradicts the assumption $a \neq 0$.\\
\quad (2)-(i) Case $a \in (\gC')_0$. Express $a = a_1e_1{}' + a_2e_2 + a_3e_3{}' + a_4e_4 + a_5{e_5}' + a_6e_6 + a_7{e_7}' , a_k\in \R$. Then,  at least one element $p$ of the following elements
$$
\begin{array}{c}
    a_4e_2 - a_2e_4, \quad a_6e_4 - a_4e_6, \quad a_3{e_1}' - a_1{e_3}',
	\quad a_7{e_5}' - a_5{e_7}' 
\end{array}
$$
satisfies $pap^{-1}=-a$. Indeed, each element $p$ above satisfies $pa=-ap$, and the norm $\N(p) $ are
$$
\begin{array}{l}
   {a_4}^2+{a_2}^2, \quad {a_6}^2+{a_4}^2,\quad -{a_3}^2-{a_1}^2, 
   \quad -{a_7}^2-{a_5}^2.
\end{array}
$$
 If $\N(p)=0$ for all $p$, then we have $a_1=\cdots=a_7=0$, that is, $a=0$, which contradicts the assumption $a \neq 0$.\\
 \quad (2)-(ii) Case $a \in (\gC^C)_0$. Express $a = a_1e_1 + \cdots + a_7e_7, a_k \in C$. Then, at least one element $p$ of the following elements 
$$
\begin{array}{c}
    a_2e_1 - a_1e_2, \quad a_3e_2 - a_2e_3,  \quad a_1e_3 - a_3e_1, 
\vspace{1mm}\\
 a_4e_3 - a_3e_4, \quad a_5e_4 - a_4e_5, \quad a_6e_5 - a_5e_6, \quad a_7e_6 - a_6e_7
\end{array}
$$
satisfies $pap^{-1}=-a$. Indeed, each element $p$ above satisfies $pa=-ap$, and the norm $\N(p)$ are
$$
\begin{array}{c}
  {a_2}^2+{a_1}^2, \quad {a_3}^2+{a_2}^2,\quad {a_1}^2+{a_3}^2,
  \vspace{1mm}\\
  {a_4}^2+{a_3}^2,\quad {a_5}^2+{a_4}^2 ,\quad {a_6}^2+{a_5}^2,\quad {a_7}^2+{a_6}^2.
\end{array}
$$
If $\N(p)=0$ for all $p$, then we have $a_1=\cdots=a_7=0$, that is, $a=0$, which contradicts the assumption $a \neq 0$.
\vspace{3mm}

{\bf Lemma 2.} {\it Let $K = \gC', \gC^C$. For  any $a, b \in K_0{}^*$ such that $(a, b) = 0$ and $\N(a) + \N(b) = 0$, there exists $p \in {K_0}^\times$} ({\it which depends on $a$ and $b$}) {\it such that $\N(pap^{-1} + b) \not= 0$.}
\vspace{1mm}

{\bf Proof.} (1) Case $K = \gC^C$. Express $a = a_1e_1 + \dots +a_7e_7, b = b_1e_1 + \dots +b_7e_7, a_k, b_k \in  C$.
\vspace{1mm}
\quad (1)-(i) When $a_kb_k \not= 0$ for some $k$. Let $p = e_k$. It holds $\N(p) = 1 \not= 0$ and 
$$
\begin{array}{l}
   \N(pap^{-1} + b) = \N(a) + \N(b) + 2(pap^{-1}, b)
\qquad
\vspace{1mm}\\
    = 0 - 2(a_1b_1 + \cdots + a_{k-1}b_{k-1} - a_kb_k + a_{k+1}b_{k+1} + \cdots + \a_7b_7)
\vspace{1mm}\\
    = 4a_kb_k +2(a, b) = 4a_kb_k \not= 0.
\end{array}$$

(1)-(ii) When $a_kb_k = 0$ for all $k = 1, 2, \cdots, 7$. There exist $k, l$ such that $k \not= l$, $a_k \not= 0, b_l \not= 0$. In this case, note that $a_l=b_k = 0$. Now, let $p = e_k + e_l$. Then, $\N(p) = 2$ and $\N(pap^{-1} + b) = \N(a) + \N(b) +2(a_kb_l - (a, b)) = 0 +2(a_kb_l - 0) = 2a_kb_l \not= 0$.
\vspace{1mm}

(2) Case $K = {\gC}'$. Express $a = a_1e_1{}' + a_2e_2 + a_3e_3{}'  + a_4e_4+ a_5e_5{}' +a_6e_6 + a_7{e_7}', b = b_1e_1{}' + b_2e_2 + b_3e_3{}'  + b_4e_4+ b_5e_5{}' +b_6e_6 + b_7{e_7}', a_k, b_k \in \R$. 
\vspace{1mm}

(2)-(i) When $a_kb_k \not= 0$ for some $k$. If $k \in \{2, 4, 6\}$ (resp. $k \in \{1,3,5,7\}$), then let $p = e_k$ (resp. $e_k{}'$). Then $\N(p) = 1$ (resp. $ -1$) $\not= 0$ and $\N(pap^{-1} + b) = 4a_kb_k $ (resp. $-4a_kb_k$) $\not= 0$ in the same way as 1-(i).
\vspace{1mm}

(2)-(ii) When $a_kb_k = 0$ for all $k = 1, 2, \cdots, 7$. There exist $k, l \in \{2, 4, 6\}$ such that $k \not= l$, $a_k \not= 0$, $b_l \not= 0$. In this case, note that $a_l=b_k = 0$. Now, let $p = e_k + e_l$. Then, $\N(p) = 2$ and $\N(pap^{-1} + b) = \N(a) + \N(b) +2(a_kb_l - (a, b)) = 0 +2(a_kb_l - 0) = 2 a_kb_l \not= 0$.
\vspace{2mm}

{\bf Theorem 3.} (1) {\it For any $a, b \in (\gC)_0{}^*$ such that $\N(a) = \N(b)$, there exists $p\in (\gC)_0{}^*$ such that $pap^{-1}=b$.} 
\vspace{1mm}

(2) {\it Let $K = \gC', \gC^C$. For any $a, b \in K_0{}^*$ such that $\N(a) = \N(b)$, there exist $p, q \in K_0{}^\times$ such that $q(pap^{-1})q^{-1}=b$.}
\vspace{3mm}

{\bf Proof.} 
(1)-(i) Case $\N(a + b) \not= 0$. Let $p = a + b$, then we have
\begin{eqnarray*}
   pap^{-1} \!\!\!&=&\!\!\! (a + b)a(a + b)^{-1}
                          = (a^2 + ba)(a + b)^{-1}
\vspace{1mm}\\
            \!\!\!&=&\!\!\! (b^2 + ba)(a + b)^{-1}
            = b(b + a)(a + b)^{-1}= b.
\end{eqnarray*}
\quad (1)-(ii)Case $\N(a + b) = 0$. This implies $b=-a$. Then, there exists $p \in (\gC)_0{}^*$ such that $pap^{-1}=-a=b$, by Lemma 1(1). 

(2)-(i) Case $\N(a + b) \not= 0$. Let $p = a + b$. Then we have $pap^{-1} =b$ in the same way as (1)-(i). 

(2)-(ii) Case $\N(a - b) \not= 0$. For given $b$, there exists $q \in {K_0}^\times$ such that $qbq^{-1}=-b$ (Lemma 1.(2)). Now, let $p = a - b$, then we have
$$
\begin{array}{l}
   q(pap^{-1})q^{-1} = q((a - b)a(a - b)^{-1})q^{-1}= q((a^2-ba)(a-b)^{-1})q^{-1}
\vspace{1mm}\\
\qquad \quad \quad \quad
 \,\,\;= q((b^2-ba)(a-b)^{-1})q^{-1} = q(b(b-a)(a-b)^{-1})q^{-1} 
\vspace{1mm}\\
\qquad \quad \quad \quad
 \,\, \;= q(-b)q^{-1}=b.
\end{array}$$

(2)-(iii) Case $\N(a + b) = \N(a - b) = 0$. This implies $(a, b) = 0$ and $\N(a) + \N(b) = 0$. Then we get $\N(a)=\N(b)=0$, that is, $a^2=b^2=0$.
For $a, b\in {K_0}^*$, we can choose $p \in {K_0}^\times$ such that 
$\N(pap^{-1} + b) \not= 0$, by Lemma 2. Let $q = pap^{-1} + b$. Then, noting that $(pap^{-1})^2=-\N(pap^{-1})=0$, we have
\begin{align}
q(pap^{-1})q^{-1} &= (pap^{-1} + b)(pap^{-1})(pap^{-1} + b)^{-1} \notag\\
&=((pap^{-1})^2 + b(pap^{-1}))(pap^{-1} + b)^{-1} \notag\\
&= (b^2 + b(pap^{-1}))(pap^{-1} + b)^{-1} \notag\\
&= (b(b + pap^{-1}))(pap^{-1} + b)^{-1}=b. \notag
\end{align}

{\bf Remark}.  In the split Cayley algebra ${\gC}'$ (resp. the complex Cayley algebra ${\gC}^C$), there exist two elements $a,b \in ({\gC}')_0{}^*$ (resp. $({\gC}^C)_0{}^*$) with the same norm $\N(a)=\N(b)$ such that $b$ can not be expressed by $pap^{-1}$ for any $p \in {\gC}'$ (resp. ${\gC}^C$).

For example, in ${\gC}'$, let $a=4e_1{}'+5e_2+3e_3{}'-5e_4+4e_5{}'+3e_7{}'$, $b=3e_2+4e_6+5e_7{}'$. Then, $\N(a)=\N(b)=0$. 
If $p=x_0+x_1e_1{}'+x_2e_2+ x_3e_3{}'+x_4{e_4}+ x_5e_5{}' +x_6e_6+ x_7{e_7}{}'\in {\gC}'$ satisfies $pa=bp$, then 
\[
\left\{
\begin{array}{@{\,}r@{\,}l@{\,}l@{\,}l@{\,}l@{\,}l@{\,}l@{\,}l@{\;}l}
{}&+4x_1 &-2 x_2 &+ 3x_ 3 &+ 5x_ 4 &+ 4x_ 5  &+ 4x_ 6  &-2x_7 &=0 \\
4x_0 &{} &+3x_ 2  &-8x_ 3  &+4x_4   &+ 5x_5   &+8x_ 6  &-4x_7 &=0 \\
2x_0  &+3x_1 &{} &-4x_3  & -4x_4   &-8x_ 5  &-5x_6  & + 4x_7 &=0 \\
3x_0  &+8x_1  &-4x_ 2  &{} &+8x_4  &+4x_ 5 &-4x_ 6 &+ 5x_7 &=0 \\
-5x_ 0 &+4x_1  &+ 4x_2  &+8x_3  &{} &-4x_5  &-8x_6&-3x_7 &=0 \\
4x_0 &-5x_ 1 &-8 x_2 &-4x_ 3 &-4x_4  &{} &+ 3x_6 &+8x_7 &=0 \\
-4x_0  &+8x_1  &+5x_2  &-4x_3 &+8x_4  &+3x_5 &{} &-4x_7 &=0 \\
-2 x_0 &+ 4x_1  &+ 4x_2  &-5x_3  &-3x_4  &-8x_5  &-4x_6 &{} &=0

\end{array}.
\right.
\]
Hence,  $p$ must be
\vspace{2mm}

$p=s(104e_1{}'+40e_2+3e_3{}'-165e_4+132e_5{}'+24e_7{}')$\\
\hspace{5cm}$+t(-46e_1{}'-8e_2+3e_3{}'+75e_4-60e_5{}'+6e_6)$
\vspace{2mm}

\noindent with arbitrary parameters, $s, t \in \R$. However, for any $s, t \in \R$, it holds $\N(p)=0$.
Therefore, there does not exist $p \in {\gC}' $ such that $pap^{-1}=b$.

Next, in ${\gC}^C$, let $a=4ie_1+5e_2+3ie_3-5e_4+4ie_5+3ie_7$, $b=3e_2+4e_6+5ie_7$. Then $\N(a)=\N(b)=0$. 
If $p=x_0+x_1e_1+x_2e_2+ x_3e_3+x_4{e_4}+ x_5{e_5} +x_6{e_6}+ x_7{e_7}\in {\gC}^C$ satisfies $pa=bp$, then 
\[
\left\{
\begin{array}{@{\,}r@{\,}l@{\,}l@{\,}l@{\,}l@{\,}l@{\,}l@{\,}l@{\;}l}
{}&-4ix_1 &-2 x_2 &-3ix_ 3 &+ 5x_ 4 &- 4ix_ 5  &+ 4x_ 6  &+2ix_7 &=0 \\
4ix_0 &{} &+3ix_ 2  &-8x_ 3  &+4ix_4   &+ 5x_5   &+8ix_ 6  &-4x_7 &=0 \\
2x_0  &-3ix_1 &{} &+4ix_3  & -4x_4   &+8ix_ 5  &-5x_6  & - 4ix_7 &=0 \\
3ix_0  &+8x_1  &-4ix_ 2  &{} &-8ix_4  &+4x_ 5 &-4ix_ 6 &+ 5x_7 &=0 \\
-5x_ 0 &-4ix_1  &+ 4x_2  &-8ix_3  &{} &+4ix_5  &-8x_6&+3ix_7 &=0 \\
4ix_0 &-5x_ 1 &-8i x_2 &-4x_ 3 &-4ix_4  &{} &+ 3ix_6 &+8x_7 &=0 \\
-4x_0  &-8ix_1  &+5x_2  &+4ix_3 &+8x_4  &-3ix_5 &{} &+4ix_7 &=0 \\
-2 ix_0 &+ 4x_1  &+ 4ix_2  &-5x_3  &-3ix_4  &-8x_5  &-4ix_6 &{} &=0
\end{array}.
\right.
\]
Hence, $p$ must be
\vspace{2mm}

$p=s(104e_1-40ie_2+3e_3+165ie_4+132e_5+24e_7)$\\
\hspace{5cm}$+t(-46ie_1-8e_2+3ie_3+75e_4-60ie_5+6e_6)$
\vspace{2mm}

\noindent with arbitrary parameters, $s, t \in C$. However, for any $s, t \in C$, it holds $\N(p)=0$.
Therefore, there does not exist $p \in {\gC}' $ such that $pap^{-1}=b$. 
\vspace{3mm}

\begin{flushleft}
{\large{\bf 2. Cases of quaternion algebras}}
\end{flushleft}

\vspace{3mm}

{\bf Lemma 4.} (1) {\it  For any $a \in (\H)_0$, there exists $p \in (\H)_0{}^* $ such that $ pap^{-1}=-a$.}\\
\quad (2) {\it  Let $K = \H'$, $\H^C$. For any $a \in K_0$, there exists $p \in K_0{}^\times $ such that $ pap^{-1}=-a$.}
\vspace{2mm}

 {\bf Proof.} Since Cayley algebras $\gC$, $\gC'$ and $\gC^C$ naturally contain quaternion algebras $\H$, $\H'$ and $\H^C$ respetively, this lemma has already been shown by Lemma 1.
\vspace{3mm}

{\bf Lemma 5.} {\it Let $K = \H', \H^C$. For any $a, b \in K_0{}^*$ such that $(a, b) = 0$ and $\N(a) +\N(b) = 0$, there exists $p \in {K_0}^\times$} ({\it which depends on $a$ and $b$}){\it \, such that $\N(pap^{-1} + b) \not= 0$.}\\
\vspace{-2mm}

{\bf Proof.} Since Cayley algebras $\gC'$ and $\gC^C$ naturally contain quaternion algebras $\H'$ and $\H^C$ respetively, this lemma has already been shown by Lemma 2.
\vspace{3mm}

{\bf Theorem 6.} 
(1) {\it For any $ a, b \in (\H)_0{}^*$ such that $\N(a) = \N(b)$, there exists $p \in (\H)_0{}^*$ such that $pap^{-1}=b$.}
\vspace{1mm}

(2) {\it Let $K=\H',\,\H^C$.  For any $ a, b \in K_0{}^*$ such that $\N(a) = \N(b)$, there exists $p \in K^\times$ such that $pap^{-1}=b $.}
\vspace{3mm}

{\bf Proof.} Since Cayley algebras $\gC$, $\gC'$ and $\gC^C$ naturally contain quaternion algebras $\H$, $\H'$ and $\H^C$ respetively, this is the particular case of Theorem 3. In the case of (2), since the associativity is valid in $\H'$ and $\H^C$, Theorem 3.(2) implies (2).

\end{document}